\documentclass[12pt, letterpaper]{article}
\usepackage{amsmath, amsfonts, amsthm, amssymb, anysize, graphicx, fullpage, setspace, hyperref, thmtools, mathrsfs}
\usepackage[authoryear, round]{natbib}
\usepackage{enumerate}
\usepackage{color}

\declaretheorem[style=definition,qed=$\square$,numberwithin=section]{definition}
\declaretheorem[style=definition,qed=$\square$,sibling=definition]{example}
\declaretheorem[style=definition,qed=$\square$,sibling=definition]{theorem}

\theoremstyle{definition}

\theoremstyle{remark}

\title{Central Limit Theorems and Large Deviations for Additive Functionals of Reflecting Diffusion Processes}
\author{Peter W. Glynn\thanks{Department of Management Science and Engineering, Stanford University. Stanford, CA, 94305. Email Address: glynn@stanford.edu.}\hspace{1 cm}
Rob J. Wang\thanks{Department of Management Science and Engineering, Stanford University. Stanford, CA, 94305. Email Address: robjwang@stanford.edu.}}
\date{October 27, 2013}
\begin{document}
\maketitle
\renewcommand{\baselinestretch}{1}
\numberwithin{equation}{section}

  \begin{abstract}
  This paper develops central limit theorems (CLT's) and large deviations results for additive functionals associated with reflecting diffusions in which the functional may include a term
	associated with the cumulative amount of boundary reflection that has occurred. Extending the known central limit and large deviations theory for Markov processes to include additive
	functionals that incorporate boundary reflection is important in many applications settings in which reflecting diffusions arise, including queueing theory and economics. In particular,
	the paper establishes the partial differential equations that must be solved in order to explicitly compute the mean and variance for the CLT, as well as the associated rate function
	for the large deviations principle.
  \end{abstract}

  \section{Introduction}
  	
	Reflecting diffusion processes arise as approximations to stochastic models associated with a wide variety of different applications domains, including communications networks,
	manufacturing systems, call centers, finance, and the study of transport phenomena (see, for example, \citet{Chen}, \citet{Harrison}, and \citet{Costantini}). If $X = (X(t):\, t \geq 0)$ 
	is the reflecting diffusion, it is often of interest
	to study the distribution of an additive functional of the form
	\[A(t) \stackrel{\Delta}{=} \int_0^t f(X(s))ds + \Lambda(t),\]
	where $f$ is a real-valued function defined on the domain of $X$, and $\Lambda = (\Lambda(t):\,t \geq 0)$ is a process (related to the boundary reflection) that increases only when
	$X$ is on the boundary of its domain. In many applications settings, the boundary process $\Lambda$ is a key quantity, as it can correspond to the cumulative number of customers lost
	in a finite buffer queue, the cumulative amount of cash injected into a firm, and other key performance measures depending on the specific application.
	
	Given such an additive functional $A = (A(t):\,t \geq 0)$, a number of limit theorems can be obtained in the setting of a positive recurrent process $X$.
	
  \vspace{0.5 cm}
  \noindent
  \textit{The Strong Law}: Compute the constant $\alpha$ such that 
  \begin{equation}\label{strong law}
  \frac{A(t)}{t} \stackrel{a.s.}{\rightarrow} \alpha
  \end{equation}
  as $t \rightarrow \infty$. In the presence of \eqref{strong law}, we can approximate $A(t)$ via 
  \begin{equation}\label{SLL approximation}
  A(t) \stackrel{\mathcal{D}}{\approx} \alpha t,
  \end{equation}
	where $\stackrel{\mathcal{D}}{\approx}$ means ``has approximately the same distribution as'' (and no other rigorous meaning, other than that supplied by
	\eqref{strong law} itself.)
	
	\vspace{0.5 cm}
	\noindent
	\textit{The Central Limit Theorem}: Compute the constants $\alpha$ and $\eta$ such that 
  \begin{equation}\label{CLT}
  t^{1/2} \left(\frac{A(t)}{t} - \alpha \right) \Rightarrow \eta N(0,\,1)
  \end{equation}
  as $t \rightarrow \infty$, where $\Rightarrow$ denotes convergence in distribution and $N(0,\,1)$ is a normal random variable (rv) with mean $0$ and unit variance.
	When \eqref{CLT} holds, we may improve the approximation \eqref{SLL approximation} to
  \begin{equation}\label{CLT approximation}
  A(t) \stackrel{\mathcal{D}}{\approx} \alpha t + \eta \sqrt{t} N(0,\,1)
  \end{equation}
	for large $t$, thereby providing a description of the distribution of $A(t)$ at scales of order $t^{1/2}$ from $\alpha t$.
	
	\vspace{0.5 cm}
	\noindent
	\textit{Large Deviations}: Compute the rate function $(I(x):\, x \in \mathbb{R})$ for which
  \begin{equation}\label{large deviations}
  \frac{1}{t} \log P(A(t) \in t \Gamma) \rightarrow -\inf_{x \in \Gamma} I(x)
  \end{equation}
  as $t \rightarrow \infty$, for subsets $\Gamma$ that are suitably chosen. Given the limit theorem \eqref{large deviations}, this suggests the (crude) approximation
  \begin{equation}\label{LD approximation}
  P(A(t) \in \Gamma) \approx \exp\left( -t \inf_{y \in \Gamma} I(y/t) \right)
  \end{equation}
	for large $t$; the approximation \eqref{LD approximation} is particularly suitable for subsets $\Gamma$ that are ``rare'' in the sense that they are more than
	order $\sqrt{t}$ from $\alpha t$.
	
	\vspace{0.5 cm}
	The main contribution of this paper concerns the computation of the quantities $\alpha$, $\eta$, and $I(\cdot)$, when $A$ is an additive functional for a reflecting diffusion
	that incorporates the boundary contribution $\Lambda$. To give a sense of the new issues that arise in this setting, observe that when $\Lambda(t) \equiv 0$ for $t \geq 0$,
	then $\alpha$ can be easily computed from the stationary distribution $\pi$ of $X$ via
	\[\alpha = \int_S f(x) \pi(dx),\]
	where $S$ is the domain of $X$. However, when $\Lambda$ is non-zero, this approach to computing $\alpha$ does not easily extend. The key to building a suitable computational theory
	for reflecting diffusions is to systematically exploit the martingale ideas that (implicitly) underly the corresponding calculations for Markov processes without boundaries; see, for example,    
  \citet{Bhattacharyya} for a discussion in the central limit setting. In the one-dimensional context, a (more laborious) approach based on the theory of regenerative processes can also be used;
	see \citet{Williams} for such a calculation in the setting of Brownian motion. In the course of our development of the appropriate martingale ideas, we will recover the existing theory
	for non-reflecting diffusions as a special case.
		
	
  The paper is organized as follows. In Section 2, we show how one can apply stochastic calculus and martingale ideas to derive partial differential equations from which the
	central limit and law of large numbers behavior for additive functionals involving boundary terms can be computed. Section 3 develops the corresponding large deviations theory
	for such additive functionals. Finally, Sections 4 and 5 illustrate the ideas in the context of one-dimensional reflecting diffusions.

  \section{Laws of Large Numbers and Central Limit Theorems}

  Let $S^o$ be a connected open set in $\mathbb{R}^d$, with $S$ and $\partial S$ denoting its closure and boundary, respectively. We assume that there exists a vector
	field $\gamma: \, \partial S \rightarrow \mathbb{R}^d$ satisfying
  \[\left\langle \gamma(x),\,n(x)\right\rangle > 0\]
  for $x \in \partial S$, where $n(x)$ is the unit inward normal to $\partial S$ at $x$ (assumed to exist). Accordingly, $\gamma(x)$ is always ``pointing'' into the interior of $S$.
	Given functions $\mu:\,S \rightarrow \mathbb{R}^d$ and $\sigma:\, S \rightarrow \mathbb{R}^{d \times d}$, we assume the existence, for each $x_0 \in S$, of a pair of continuous
	processes $X = (X(t):\,t\geq 0)$ and $k = (k(t):\,t \geq 0)$ (with $k$ of bounded variation) for which
  \begin{equation}\label{SDE}
  X(t) = x_0 + \int_0^t \mu(X(s))ds + \int_0^t \sigma(X(s))dB(s) + k(t),
  \end{equation}
	\[X(t)\in S,\]
  \[|k|(t) = \int_0^t I(X(r) \in \partial S)d|k|(r),\]
	and
	\[k(t) = \int_0^t \gamma(X(s))d|k|(s),\]
	where $B = (B(t):\,t\geq 0)$ is a standard $\mathbb{R}^d$-valued Brownian motion, and $|k|(t)$ is the (scalar) total variation of $k$ over $[0,\,t]$; sufficient conditions surrounding
	existence of such processes can be found in \citet{Lions}. Note that our formulation permits the direction of reflection to be oblique. Regarding the structure of the boundary process
	$\Lambda$, we assume that it takes the form
  \[\Lambda(t) = \int_0^t r(X(s)) d|k|(s),\]
  for a given function $r:\,S \rightarrow \mathbb{R}$.
	
	We expect laws of large numbers and central limit theorems to hold with the conventional normalizations only when $X$ is a positive recurrent Markov process. In view of this, we assume:
	
	\vspace{0.5 cm}
	\noindent
  A1: $X$ is a Markov process with a stationary distribution $\pi$ that is recurrent in the sense of Harris (by Harris recurrence, we mean that there exists a non-trivial
	$\sigma$-finite measure $\phi$ on $S$ for which whenever $\phi(B) > 0$, $\int_0^\infty I(X(s) \in B)ds = \infty$ $P_x$ a.s. for each $x \in S$, where
  \[P_x(\cdot) \stackrel{\Delta}{=} P(\cdot\,|\,X(0)=x).)\]
  
	\noindent
	We note that Harris recurrence implies that any stationary distribution must be unique. For a discussion of methods for verification of recurrence in the setting of
	continuous-time Markov processes, see \citet{Meyn}, \citet{Meyn2}, and \citet{Down}.

  The key to developing laws of large numbers and central limit theorems for the additive functional $A$ is to find a function $u:\,S \rightarrow \mathbb{R}$ and a constant $\alpha$ for which
  \[M(t) \stackrel{\Delta}{=} u(X(t))-(A(t)-\alpha t)\]
  is a local $\mathcal{F}_t$-martingale, where $\mathcal{F}_t = \sigma(X(s):\,0\leq s\leq t).$ In order to explicitly compute $u$, it is convenient to identify a suitable partial
	differential equation satisfied by $u$ that can be used to solve for $u$. Note that if $u \in C^2(S)$, It\^o's formula ensures that
	
  \begin{eqnarray*}
  dM(t)
  &=& du(X(t)) - (f(X(t))-\alpha)dt - r(X(t))d|k|(t)\\
  &=& \nabla u(X(t))dX(t) + \frac{1}{2} \sum_{i,\,j=1}^d (\sigma \sigma^T)_{ij}(X(t)) \frac{\partial^2 u(X(t))}{\partial x_i \partial x_j}dt\\
  &\qquad-& (f(X(t))-\alpha) dt - r(X(t))d|k|(t)\\
  &=& \nabla u(X(t))(\mu(X(t))dt + \sigma(X(t))dB(t)\\
  &\qquad+& \gamma(X(t))d|k|(t)) + \frac{1}{2} \sum_{i,\,j=1}^d (\sigma \sigma^T)_{ij}(X(t)) \frac{\partial^2 u(X(t))}{\partial x_i \partial x_j}dt\\
  &\qquad-& (f(X(t))-\alpha) dt - r(X(t))d|k|(t)\\
  &=& ((\mathcal{L}u)(X(t))-(f(X(t))-\alpha))dt + (\nabla u(X(t)) \gamma(X(t))-r(X(t)))d|k|(t)\\
  &\qquad +& \nabla u(X(t))\sigma(X(t))dB(t),
  \end{eqnarray*}
  where $\nabla u(x)$ is the gradient of $u$ evaluated at $x$ (encoded as a row vector) and $\mathcal{L}$ is the elliptic differential operator
  \[\mathcal{L} = \sum_{i=1}^d \mu_i(x) \frac{\partial}{\partial x_i} + 
  \frac{1}{2} \sum_{i,\,j=1}^d (\sigma \sigma^T)_{ij}(x) \frac{\partial^2}{ \partial x_i \partial x_j}.\]
  The process $M$ can be guaranteed to be a local martingale if we require that $u$ and $\alpha$ satisfy
  \begin{gather}\label{PDE}
  (\mathcal{L} u)(x) = f(x) - \alpha,\;x\in S\\
  \nabla u(x) \gamma(x) = r(x),\; x \in \partial S,\notag
  \end{gather}
  since this choice implies that
  \[dM(t) = \nabla u(X(t)) \sigma(X(t))dB(t).\]
  (We use here the fact that $|k|(t)$ increases only when $X(t) \in \partial S$.) Accordingly, the quadratic variation of $M$ is given by
  \begin{eqnarray*}
  [M,\,M](t)
  &=& \int_0^t \nabla u(X(s))\sigma(X(s))\sigma^T(X(s)) \nabla u(X(s))^T ds\\
  &\stackrel{\Delta}{=}& \int_0^t \nu(X(s))ds.
  \end{eqnarray*}
  Since $\nu$ is nonnegative and $X$ is positive Harris recurrent, it follows that
  \[\frac{1}{t} \int_0^t \nu(X(s))ds \rightarrow \int_S \nu(y) \pi(dy)\;\; P_x\;a.s.\]
  as $t \rightarrow \infty$, for each $x \in S$. Set
  \begin{eqnarray*}
  \eta^2
  &=& \int_S \nu(y) \pi(dy)\\
  &=& \int_S \nabla u(y)\sigma(y)\sigma(y)^T \nabla u(y) \pi(dy),
  \end{eqnarray*}
  and assume $\eta^2 < \infty$. As a consequence of the path continuity of $M$, the martingale central limit theorem then implies that for each $x \in S$, 
  \begin{equation}\label{CLT for M}
  t^{-1/2} M(t) \Rightarrow \eta N(0,\,1)
  \end{equation}
  as $t \rightarrow \infty$ under $P_x$ (see, for example, \citet{Ethier}). In other words, 
  \[t^{-1/2} (u(X(t))-(A(t)-\alpha t)) \Rightarrow \eta N(0,\,1)\]
  as $t \rightarrow \infty$ under $P_x$.

  Let $P_\pi(\cdot) = \int_S P_x(\cdot)\pi(dx)$, and observe that $X$ is stationary under $P_\pi$. Thus, $u(X(t)) \stackrel{\mathcal{D}}{=} u(X(0))$ for $t \geq 0$ under
	$P_\pi$ (where $\stackrel{\mathcal{D}}{=}$ denotes equality in distribution), so that
  \begin{equation}\label{vanishing}
  t^{-1/2} u(X(t)) \Rightarrow 0
  \end{equation}
  as $t \rightarrow \infty$ under $P_\pi$. It follows that
  \[\frac{1}{t} A(t) \Rightarrow \alpha\]
  as $t \rightarrow \infty$ under $P_\pi$. Let $E_\pi (\cdot)$ be the expectation operator associated with $P_\pi$. If $f$ and $r$ are nonnegative, the Harris recurrence implies that
  \[\frac{1}{t} A(t) \rightarrow E_\pi A(1)\;\;P_x\;a.s.\]
  as $t \rightarrow \infty$, for each $x \in S$. Hence, $E_\pi A(1) = \alpha$, so that
  \[\frac{1}{t} A(t) \rightarrow \alpha\;\;P_x\;a.s.\]
  as $t \rightarrow \infty$, for each $x \in S$. This establishes the desired strong law of large numbers for the additive functional $A$.

  Turning now to the central limit theorem, \eqref{CLT for M} and \eqref{vanishing} together imply that
  \[t^{1/2} \left( \frac{A(t)}{t} - \alpha \right) \Rightarrow \eta N(0,\,1)\]
  as $t \rightarrow \infty$ under $P_\pi$. Since $X$ is a one-dependent regenerative process (i.e. there exists a strictly increasing sequence
	$(T_n:\,n \geq -1)$ of random times for which $T_{-1}=0$ and $(W_j:\,j \geq 0)$ is a one-dependent sequence of random elements and $(W_j:\,j \geq 1)$
	is an identically distributed sequence, where
	\[W_j \stackrel{\Delta}{=} (X(T_{j-1} +s): 0 \leq s < T_j - T_{j-1});\]
	see \citet{Glynn} and \citet{Sigman} for further details), evidently
  \[t^{1/2} \left( \frac{A(t)}{t} - \alpha \right) \Rightarrow \eta N(0,\,1)\]
  as $t \rightarrow \infty$ under $P_x$, for each $x \in S$. We summarize this discussion with the following theorem.

  \begin{theorem}\label{main theorem}
  Assume A1 and that $f$ and $r$ are nonnegative. If there exists $u \in C^2(S)$ and $\alpha \in \mathbb{R}$ that satisfy
  \begin{gather*}
  (\mathcal{L} u)(x) = f(x) - \alpha,\;x\in S\\
  \nabla u(x) \gamma(x) = r(x),\; x \in \partial S,
  \end{gather*}
  with
  \[\eta^2 = \int_S \nabla u(y)\sigma(y)\sigma(y)^T \nabla u(y) \pi(dy) < \infty,\]
  then, for each $x \in S$,
  \[\frac{1}{t} A(t) \rightarrow \alpha\;\;P_x\;a.s.\]
  and
  \[t^{1/2} \left( \frac{A(t)}{t} - \alpha \right) \Rightarrow \eta N(0,\,1)\]
  as $t \rightarrow \infty$, under $P_x$.
  \end{theorem}
  The function $u$ satisfying \eqref{PDE} is said to be a solution of the \textit{generalized Poisson equation} corresponding to the pair $(f,\,r)$.

  \section{Large Deviations for the Additive Functional \texorpdfstring{$A$}{TEXT}}

  The key to developing a suitable large deviations theory for $A$ is again based on construction of an appropriate martingale. Here, we propose a one-parameter family of martingales of the form
  \[M(\theta,\,t) = \exp(\theta A(t)-\psi(\theta)t)h_\theta(X(t))\]
  for $\theta$ lying in some open interval containing the origin, where $\psi(\theta)$ and $h_\theta$ are chosen
  appropriately. As in Section 2, we use stochastic calculus to derive a corresponding PDE from which one can potentially compute $\psi(\theta)$ and $h_\theta$ analytically.
	In particular, if $h_\theta \in C^2(S)$, It\^o's formula yields
	
  \begin{eqnarray*}
  dM(\theta,\,t)
  &=& d(\exp(\theta A(t)-\psi(\theta)t))h_\theta(X(t))\\
  &\qquad+& \exp(\theta A(t)-\psi(\theta)t)dh_\theta(X(t))\\
  &=& \exp(\theta A(t)-\psi(\theta)t) (\theta f(X(t))dt + \theta r(X(t))d|k|(t) - \psi(\theta)dt) h_\theta(X(t))\\
  &\qquad+&  \exp(\theta A(t)-\psi(\theta)t)\bigg[\nabla h_\theta(X(t))\mu(X(t))dt + \nabla h_\theta(X(t))\sigma(X(t))dB(t)\\
  &\qquad+& \nabla h_\theta(X(t)) \gamma(X(t))d|k|(t) + \frac{1}{2} \sum_{i,\,j=1}^d (\sigma \sigma^T)_{ij}(X(t)) \frac{\partial^2 h_\theta(X(t))}{\partial x_i \partial x_j} dt\bigg]\\
  &=& \exp(\theta A(t)-\psi(\theta)t) \bigg[((\mathcal{L}h_\theta)(X(t)) + (\theta f(X(t))-\psi(\theta))h_\theta(X(t)))dt\\
  &\qquad+& (\nabla h_\theta(X(t))\gamma(X(t)) + \theta r(X(t))h_\theta(X(t)))d|k|(t)\\
  &\qquad+& \nabla h_\theta(X(t))\sigma(X(t))dB(t)\bigg],
  \end{eqnarray*}
  where $\mathcal{L}$ is the differential operator defined in Section 2. If we require that $h_\theta$ and $\psi(\theta)$ satisfy
  \begin{gather}\label{LD PDE}
  (\mathcal{L} h_\theta)(x) + (\theta f(x)-\psi(\theta))h_\theta(x) = 0,\;x\in S\\
  \nabla h_\theta(x) \gamma(x) + \theta r(x)h_\theta(x) = 0,\; x \in \partial S,\notag
  \end{gather}
  then
  \[dM(\theta,\,t) = \nabla h_\theta(X(t))\sigma(X(t))dB(t),\]
  and $M(\theta,\,t):\,t \geq 0)$ is consequently a local $\mathcal{F}_t$-martingale. (Again, we use here the fact that $|k|$ increases only when $X$ is on the boundary of $S$.)
	Note that \eqref{LD PDE} takes the form of an eigenvalue problem involving the operator $\mathcal{L} + \theta f \mathcal{I}$, where $\mathcal{I}$ is the identity operator for which
	$\mathcal{I} u = u$. In this eigenvalue formulation, $\psi(\theta)$ is the eigenvalue and $h_\theta$ the corresponding eigenfunction. Since $\mathcal{L} + \theta f \mathcal{I}$ is
	expected to have multiple eigenvalues, \eqref{LD PDE} cannot be expected to uniquely determine $\psi(\theta)$ and $h_\theta$. In order to ensure uniqueness, we now add the requirement
	that $h_\theta$ be positive.

  Let $(T_n:\,n \geq 0)$ be the localizing sequence of stopping times associated with the local martingale $(M(\theta,\,t):\,t\geq 0)$, so that
  \begin{equation}\label{localization}
  E_x \exp(\theta A(t \wedge T_n) - \psi(\theta)(t \wedge T_n)) h_\theta(X(t \wedge T_n)) = h_\theta(x)
  \end{equation}
  for $x \in S$, where $E_x(\cdot)$ is the expectation operator associated with $P_x(\cdot)$ and $a \wedge b 
  \stackrel{\Delta}{=} \min(a,\,b)$ for $a,\,b \in \mathbb{R}$.

  Suppose that $S$ is compact, so that $h_\theta$ is then bounded above and below by positive constants (on account of the positivity of $h_\theta$ and the fact that $h_\theta \in C^2(S)$).
	If $f$ and $r$ are nonnegative (as in Section 2), it follows that for $\theta \leq 0$,
  \[\exp(\theta A(t \wedge T_n) - \psi(\theta)(t \wedge T_n)) h_\theta(X(t \wedge T_n))\]
  is a bounded sequence of rv's, and thus the Bounded Convergence Theorem implies that
  \begin{equation}\label{negative theta}
  E_x \exp(\theta A(t) - \psi(\theta)t) h_\theta(X(t)) = h_\theta(x)
  \end{equation}
  for $\theta \leq 0$, and $x \in S$.

  On the other hand, if $\theta > 0$, the positivity of $h_\theta$ and Fatou's lemma imply that
  \[E_x \exp(\theta A(t) - \psi(\theta)t) h_\theta(X(t)) \leq h_\theta(x)\]
  for $x \in S$, from which we may obtain the upper bound
  \[E_x \exp(\theta A(t)) \leq e^{\psi(\theta)t} \frac{h_\theta(x)}{\inf_{y \in S} h_\theta(y)},\]
  and hence $\exp(\theta A(t))$ is $P_x$-integrable. Since $f$ and $r$ are nonnegative and $\theta > 0$,
  $\theta A(t \wedge T_n) \leq \theta A(t)$, so
  \[\exp(\theta A(t \wedge T_n) - \psi(\theta)(t \wedge T_n)) h_\theta(X(t \wedge T_n))
  \leq \exp(\theta A(t) + |\psi(\theta)|(t)) \sup_{y \in S} h_\theta(y).\]
  The Dominated Convergence Theorem, as applied to \eqref{localization}, then yields the conclusion that
  \begin{equation}\label{positive theta}
  E_x \exp(\theta A(t) - \psi(\theta)t) h_\theta(X(t)) = h_\theta(x)
  \end{equation}
  for $x \in S$. Since
  \[e^{\psi(\theta)t} \frac{h_\theta(x)}{\sup_{y \in S} h_\theta(y)} \leq E_x \exp(\theta A(t)) \leq
  e^{\psi(\theta)t} \frac{h_\theta(x)}{\inf_{y \in S} h_\theta(y)},\]
  it follows that
  \[\frac{1}{t} \log E_x \exp(\theta A(t)) \rightarrow \psi(\theta)\]
  as $t \rightarrow \infty$, proving the following theorem.

  \begin{theorem}\label{LD Theorem}
  Assume that $S$ is compact and that $f$ and $r$ are nonnegative. If there exists a positive function $h_\theta \in C^2(S)$ and $\psi(\theta) \in \mathbb{R}$ that satisfy
  \begin{gather*}
  (\mathcal{L} h_\theta)(x) + (\theta f(x)-\psi(\theta))h_\theta(x) = 0,\;x\in S\\
  \nabla h_\theta(x) \gamma(x) + \theta r(x)h_\theta(x) = 0,\; x \in \partial S,
  \end{gather*}
  then
  \[\frac{1}{t} \log E_x \exp(\theta A(t)) \rightarrow \psi(\theta)\]
  as $t \rightarrow \infty$.
  \end{theorem}
  The G\"{a}rtner-Ellis Theorem (see, for example, p.45 of \citet{Dembo}) then provides technical conditions under which
  \[\frac{1}{t} \log P_x(A(t) \in t \Gamma) \rightarrow -\inf_{y \in \Gamma} I(y)\]
  as $t \rightarrow \infty$, where
  \[I(y) = \sup_{\theta \in \mathbb{R}} [\theta y - \psi(\theta)].\]
  In particular, if $\Gamma = (z,\,\infty)$, then
  \[\frac{1}{t} \log P_x(A(t) \geq tz) \rightarrow -(\theta_z z - \psi(\theta_z)),\]
  provided that $\psi(\cdot)$ is differentiable and strictly convex in a neighborhood of a point $\theta_z$ satisfying $\psi'(\theta_z) = z$.
	See p.15-16 of \citet{Bucklew} for a related argument.

  \section{CLT's for One-dimensional Reflecting Diffusions}

  We now illustrate these ideas in the setting of one-dimensional diffusions. In this context, we can compute the solution of the generalized Poisson equation corresponding to $(f,\,r)$
	fairly explicitly.

  We start with the case where there are two reflecting barriers, at $0$ and $b$, so that $S = [0,\,b]$.
  Then, $X = (X(t):\,t\geq 0)$ satisfies the stochastic differential equation (SDE)
  \begin{eqnarray*}
  dX(t) 
  &=& \mu(X(t))dt + \sigma(X(t))dB(t) + dL(t) - dU(t)\\
  &=& \mu(X(t))dt + \sigma(X(t))dB(t) + dk(t),
  \end{eqnarray*}
  with $\gamma(0) = 1$ and $\gamma(b) = -1$; the processes $L$ and $U$ increase only when $X$ visits the lower and
  upper boundaries at $0$ and $b$, respectively. We consider here the additive functional
  \[A(t) = \int_0^t f(X(s))ds + r_0 L(t) + r_b U(t),\]
  where $f:[0,\,b]\rightarrow \mathbb{R}$ is assumed to be bounded. In this setting, Theorem \ref{main theorem}
  leads to consideration of the ordinary differential equation (ODE)
  \begin{equation}\label{one-dimensional ODE}
  \mu(x) u'(x) + \frac{\sigma^2(x)}{2} u''(x) = f(x) - \alpha,
  \end{equation}
  \begin{equation}
  u'(0) = r_0,
  \end{equation}
  \begin{equation}
  u'(b) = -r_b.
  \end{equation}

  Hence, if $\mu(\cdot)$ and $\sigma^2(\cdot)$ are continuous and $\sigma^2(\cdot)$ positive, \eqref{one-dimensional ODE} can be re-written via the method of integrating factors
	(see, for example, \citet{Karlin}) as
  \[\frac{d}{dx} \left( \exp\left( \int_0^x \frac{2\mu(y)}{\sigma^2(y)}dy \right) u'(x)\right)
  = \frac{2(f(x)-\alpha)}{\sigma^2(x)} \exp\left( \int_0^x \frac{2\mu(y)}{\sigma^2(y)}dy \right),\]
  from which we conclude that
  \begin{eqnarray}\label{u prime}
  u'(x)
  &=& \left( u'(0) + \int_0^x \frac{2(f(y)-\alpha)}{\sigma^2(y)} \exp\left( \int_0^y \frac{2\mu(z)}{\sigma^2(z)}dz \right)dy \right)\\
  &\quad \cdot& \exp\left( -\int_0^x \frac{2\mu(y)}{\sigma^2(y)}dy \right).\notag
  \end{eqnarray}
  But $u'(0) = r_0$ and $u'(b) = -r_b$, and thus
  \begin{eqnarray*}
  -r_b
  &=& \left( r_0 + \int_0^b \frac{2(f(y)-\alpha)}{\sigma^2(y)} \exp\left( \int_0^y \frac{2\mu(z)}{\sigma^2(z)}dz \right)dy \right)\\
  &\quad \cdot& \exp\left( -\int_0^b \frac{2\mu(y)}{\sigma^2(y)}dy \right).
  \end{eqnarray*}
  Hence,
  \begin{equation}\label{alpha}
  \alpha = \frac{r_0 + r_b e^{\left( \int_0^b \frac{2\mu(y)}{\sigma^2(y)}dy \right)} + \int_0^b
  \frac{2f(y)}{\sigma^2(y)} e^{\left( \int_0^y \frac{2\mu(z)}{\sigma^2(z)}dz \right)}dy}{2 \int_0^b \frac{1}{\sigma^2(y)} e^{\left( \int_0^y \frac{2\mu(z)}{\sigma^2(z)}dz \right)}dy}
  \end{equation}
  By setting $r_0 = r_b = 0$, we conclude that the stationary distribution $\pi$ of $X$ must satisfy
  \begin{equation}\label{stationary distribution}
  \int_0^b \pi(dx)f(x) = \int_0^b f(x)p(x)dx,
  \end{equation}
  where
  \[p(x) = \frac{\frac{1}{\sigma^2(x)} \exp{\left( \int_0^x \frac{2\mu(z)}{\sigma^2(z)}dz \right)}}{
  \int_0^b \frac{1}{\sigma^2(y)}\exp{\left( \int_0^y \frac{2\mu(z)}{\sigma^2(z)}dz \right)}dy}.\]
  Since \eqref{stationary distribution} holds for all bounded functions $f$, it follows that $\pi(dx) = p(x)dx$, so that $\pi$ has now been computed. Furthermore, \eqref{u prime}
  establishes that
  \begin{eqnarray*}
  u'(x)
  &=& \left( r_0 + \int_0^x \frac{2(f(y)-\alpha)}{\sigma^2(y)} \exp\left( \int_0^y \frac{2\mu(z)}{\sigma^2(z)}dz \right)dy \right)\\
  &\quad \cdot& \exp\left( -\int_0^x \frac{2\mu(y)}{\sigma^2(y)}dy \right),
  \end{eqnarray*}
  where $\alpha$ is given by \eqref{alpha}. Consequently, we have explicit formulae for both $\pi$ and $u'$, from which the variance constant
  \[\eta^2 = \int_0^b u'(x)^2 \sigma^2(x)p(x)dx\]
  of Theorem \ref{main theorem} can now be calculated. We now illustrate these calculations in the context of some special cases,
  focusing our interest on the boundary processes (by setting $f = 0$).

  \begin{example}{\textit{Two-sided Reflecting Brownian Motion:}}
  Here $\mu(x) = \mu$ and $\sigma^2(x) = \sigma^2 > 0$. 
  If $\mu \neq 0$, then, upon setting $\xi = 2\mu/\sigma^2$,
  \begin{eqnarray*}
  \alpha
  &=& \frac{\mu(r_0 + r_b e^{\xi b})}{e^{\xi b} - 1}
  \end{eqnarray*}
  and
  \begin{eqnarray*}
  p(x)
  &=& \frac{\xi e^{\xi x}}{e^{\xi b}-1}.
  \end{eqnarray*}
  Also,
  \begin{eqnarray*}
  u'(x)
  &=& \left(r_0 + \int_0^x -\frac{2\alpha}{\sigma^2} e^{\frac{2\mu}{\sigma^2}y}dy\right)e^{-\frac{2\mu}{\sigma^2}x}\\
  &=& \left( \frac{r_0 + r_b}{1-e^{-\xi b}}\right)e^{-\xi x} -  \frac{r_0e^{-\xi b} + r_b }{1-e^{-\xi b}}
  \end{eqnarray*}
  and consequently
  \[u'(x)^2 = \left( \frac{r_0 + r_b}{1-e^{-\xi b}}\right)^2 e^{-2\xi x}
  -\frac{2(r_0e^{-\xi b} + r_b)(r_0 + r_b) }{(1-e^{-\xi b})^2} e^{-\xi x} + \left(\frac{r_0e^{-\xi b} + r_b }{1-e^{-\xi b}}\right)^2.\]
  Therefore,
  \begin{eqnarray*}
  \eta^2
  &=& \sigma^2 \left[ \left(\frac{r_0 + r_b}{1-e^{-\xi b}}\right)^2 e^{-\xi b}
  - \frac{(r_0e^{-\xi b} + r_b)(r_0 + r_b) }{(1-e^{-\xi b})^2} \frac{2 \xi b}{{e^{\xi b}-1}} + \left(\frac{r_0e^{-\xi b} + r_b }{1-e^{-\xi b}}\right)^2\right].
  \end{eqnarray*}

  If $\mu = 0$, then $\alpha = \frac{\sigma^2 (r_0 + r_b)}{2b}$ and $p(x) = \frac{1}{b}$. Also,
  \begin{eqnarray*}
  u'(x)
  &=& r_0 - \frac{(r_0+r_b)}{b}x
  \end{eqnarray*}
  and therefore
  \begin{eqnarray*}
  \eta^2
  &=& \sigma^2 \int_0^b \frac{\left( \frac{(r_0+r_b)}{b}x - r_0 \right)^2}{b}dx\\
  &=& \frac{\sigma^2(r_0^3 + r_b^3)}{3(r_0+r_b)}.
  \end{eqnarray*}
  \end{example}

  \begin{example}{\textit{Two-sided Reflecting Ornstein-Uhlenbeck:}}
  For this process, $\mu(x) = -a(x-c)$ and $\sigma^2(x) = \sigma^2 > 0$. We thus have
  \begin{eqnarray*}
  \alpha
  &=& \frac{r_0 + r_b e^{-\frac{a(b-c)^2-ac^2}{\sigma^2}}}{\frac{2}{\sigma^2} \int_0^b e^{-\frac{a(y-c)^2-ac^2}{\sigma^2}}dy}.
  \end{eqnarray*}
  Also,
  \begin{eqnarray*}
  u'(x)
  &=& \left( r_0 - \frac{2\alpha}{\sigma^2}\int_0^x e^{ -\int_0^y \frac{2a(z-c)}{\sigma^2}dz }dy\right)e^{\int_0^x \frac{2a(y-c)}{\sigma^2}dy}\\
  &=& r_0 e^{\frac{a(x-c)^2-ac^2}{\sigma^2}} - \frac{2\alpha}{\sigma^2} \int_0^x e^{-\frac{a(y-c)^2-a(x-c)^2}{\sigma^2}} dy
  \end{eqnarray*}
  and
  \begin{eqnarray*}
  p(x)
  &=& \frac{e^{-\int_0^x \frac{2a(z-c)}{\sigma^2}}dz}{\int_0^b e^{-\int_0^y \frac{2a(z-c)}{\sigma^2}dz}dy}\\
  &=& \sqrt{\frac{2a}{\sigma^2}} \frac{\phi\left( (x-c)\sqrt{\frac{2a}{\sigma^2}} \right)}{ \Phi\left( (b-c)\sqrt{\frac{2a}{\sigma^2}} \right) - \Phi\left( (-c)\sqrt{\frac{2a}{\sigma^2}}  \right) },
  \end{eqnarray*}
  where $\phi$ and $\Phi$ are, respectively, the density and cumulative density function (CDF) of a standard normal random variable.
  From these, one may readily compute
  \[\eta^2 = \sigma^2\int_0^b \left( r_0 e^{\frac{a(x-c)^2-ac^2}{\sigma^2}} - \frac{2\alpha}{\sigma^2} \int_0^x e^{-\frac{a(y-c)^2-a(x-c)^2}{\sigma^2}} dy \right)^2 p(x)dx\]
  numerically when the problem data are explicit.
  \end{example}

  The diffusions in our examples arise as approximations to queues in heavy traffic, in which $L(t)$ then approximates the cumulative lost service capacity of the server over $[0,\,t]$,
	while $U(t)$ describes the cumulative number of customers lost due to blocking (because of arrival to a full buffer); see \citet{Zhang} for details.

  Turning now to the setting in which only a single reflecting barrier is present (say, at the origin), $S$ then takes the form $S = [0,\,\infty)$, and the differential equation
	for $u$ takes the form
  \begin{gather*}
  \mu(x) u'(x) + \frac{\sigma^2(x)}{2}u''(x) = f(x) - \alpha,\\
  u'(0) = r_0.
  \end{gather*}
  Then $u'(x)$ is again given by \eqref{u prime}, and
  \begin{equation}\label{alpha with single barrier}
  \alpha = \frac{r_0 + \int_0^\infty \frac{2f(y)}{\sigma^2(y)} e^{\left( \int_0^y \frac{2\mu(z)}{\sigma^2(z)}dz \right)}dy}{2 \int_0^\infty \frac{1}{\sigma^2(y)}
	e^{\left( \int_0^y \frac{2\mu(z)}{\sigma^2(z)}dz \right)}dy},
  \end{equation}
  provided that the problem data are such that the integrals in \eqref{alpha with single barrier} converge and are finite. In particular, $X$ fails to have a stationary distribution if
  \[\int_0^\infty \frac{1}{\sigma^2(y)} e^{\left( \int_0^y \frac{2\mu(z)}{\sigma^2(z)}dz \right)}dy = \infty.\]

  \section{Large Deviations: One-dimensional Reflecting Diffusions}

  In this setting, we discuss the large deviations theory of Section 3, specialized to the setting of one-dimensional diffusions with reflecting barriers at $0$ and $b$.
	Theorem \ref{LD Theorem} asserts that the key ODE in this setting requires finding $\psi(\theta) \in \mathbb{R}$ and $h_\theta \in C^2[0,\,b]$ for which
  \begin{equation}\label{one-dim LD ODE}
  \mu(x) h_\theta'(x) + \frac{\sigma^2(x)}{2}h_\theta''(x) + (\theta f(x) - \psi(\theta)) h_\theta(x) = 0,\; 0 \leq x \leq b
  \end{equation}
  \begin{gather*}
  h_\theta'(0) + \theta r_0 h_\theta(0) = 0,\\
  -h_\theta'(b) + \theta r_b h_\theta(b) = 0.
  \end{gather*}
  The above differential equation \eqref{one-dim LD ODE} can be put in the form
  \begin{equation}\label{SL form}
  -\frac{d}{dx} (a(x) h'_\theta(x)) + b(x)h_\theta(x) = \lambda c(x) h_\theta(x)
  \end{equation}
  for $0 \leq x \leq b$, where $\lambda = -\psi(\theta)$ and 
  \[a(x) = \exp\left( \int_0^x \frac{2\mu(y)}{\sigma^2(y)}dy \right),\]
  \[b(x) = -\frac{2\theta f(x)}{\sigma^2(x)}  \exp\left( \int_0^x \frac{2\mu(y)}{\sigma^2(y)}dy \right),\]
  \[c(x) = \frac{2}{\sigma^2(x)}\exp\left( \int_0^x \frac{2\mu(y)}{\sigma^2(y)}dy \right).\]
  Suppose that $f$, $\mu$, and $\sigma^2$ are continuous on $[0,\,b]$, with $\sigma^2(x) > 0$ for $x \in [0,\,b]$. Because $a(\cdot)$ and $c(\cdot)$ are then positive on $[0,\,b]$,
	\eqref{SL form} takes the form of a so-called Sturm-Liouville problem. Consequently, there exist real eigenvalues $\lambda_1 < \lambda_2 < \cdots$ with $\lambda_n \rightarrow \infty$
	satisfying \eqref{SL form}, with corresponding eigenfunction solutions $v_1,\,v_2,\;\ldots$. Furthermore, the eigenfunction $v_i$ has the property that it has exactly $i-1$ roots in $[0,\,b]$;
	see, for example, \citet{Al-Gwaiz} for details on Sturm-Liouville theory. As a consequence, the eigenfunction $v_1$ is the only eigenfunction that can be taken to be positive over
	$[0,\,b]$. Thus, it follows that we should set $\psi(\theta) = -\lambda_1$ and $h_\theta = v_1$.

  We now illustrate these ideas in the setting of reflecting Brownian motion in one dimension, again focusing on the boundary process by setting $f = 0$.

  \begin{example}{\textit{Two-sided Reflecting Brownian Motion:}}
  Here $\mu(x) = \mu$ and $\sigma^2(x) = \sigma^2 > 0$. The case in which $r_0 = 0$ and $r_b = 1$ was studied in detail in \citet{Zhang}.
	In particular, consider the parameter spaces given by
	\begin{eqnarray*}
    \mathscr{R}_1 & = & \{(\theta,\,\mu,\,b):\,\theta > 0\} \\
    \mathscr{R}_2 & = & \{(\theta,\,\mu,\,b):\,\theta < 0,\, \mu(\mu+\theta \sigma^2) \leq 0\} \\
    \mathscr{R}_3 & = & \{(\theta,\,\mu,\,b):\,\theta < 0,\, \mu(\mu+\theta \sigma^2) > 0,\, b\mu(\mu+\theta \sigma^2) > -\theta \sigma^4\}\\
		\mathscr{R}_4 & = & \{(\theta,\,\mu,\,b):\,\theta < 0,\, \mu(\mu+\theta \sigma^2) > 0,\, b\mu(\mu+\theta \sigma^2) < -\theta \sigma^4\}\\
		\mathscr{B}_1 & = & \{(\theta,\,\mu,\,b):\,\theta = 0\}\\
		\mathscr{B}_2 & = & \{(\theta,\,\mu,\,b):\,\theta < 0,\, \mu(\mu+\theta \sigma^2) > 0,\, b\mu(\mu+\theta \sigma^2) = -\theta \sigma^4\}.
  \end{eqnarray*}
	The authors showed that, for $(\theta,\,\mu,\,b) \in \mathscr{R}_i$ ($i = 1,\,3$), the solutions $\psi = \psi(\theta)$ and $h_\theta(\cdot)$ to
	\begin{eqnarray*}
	(\mathcal{L}h_\theta)(x) &=& \psi(\theta) h_\theta(x)\\
	h_\theta'(0) &=& 0\\
	h_\theta'(b) &=& \theta h_\theta(b)\\
	h_\theta(0) &=& 1
	\end{eqnarray*}
	for $0 \leq x \leq b$ are given by $\psi(\theta) = \frac{\beta(\theta)^2 - \mu^2}{2\sigma^2}$ and
	\[h_\theta(x) = \frac{1}{2\beta(\theta)} e^{-\frac{\mu}{\sigma^2}x}\left[ (\beta(\theta)-\mu) e^{-\frac{\beta(\theta)}{\sigma^2} x} +  (\beta(\theta)+\mu) e^{\frac{\beta(\theta)}{\sigma^2} x}\right],\]
	where $\beta(\theta)$ is the unique root in $\mathscr{F}_i$ of the equation
	\[\frac{1}{\beta} \log \left( \frac{(\beta-\mu)(\beta+\mu+\theta \sigma^2)}{(\beta+\mu)(\beta-\mu-\theta \sigma^2)} \right) = \frac{2b}{\sigma^2},\]
	with $\mathscr{F}_1 = (|\mu| \vee |\mu+\theta \sigma^2|,\,\infty)$ and $\mathscr{F}_3 = (0,\,|\mu| \wedge |\mu+\theta \sigma^2|)$.
	For $(\theta,\,\mu,\,b) \in \mathscr{R}_i$ ($i = 2,\,4$),  the solutions are given by $\psi(\theta) = -\frac{\xi(\theta)^2 + \mu^2}{2\sigma^2}$ and
	\[h_\theta(x) = e^{-\frac{\mu}{\sigma^2}x} \left[ \cos\left( \frac{\xi(\theta)x}{\sigma^2} \right) + \frac{\mu}{\xi(\theta)} \sin\left(\frac{\xi(\theta)x}{\sigma^2}\right) \right],\]
	where $\xi(\theta)$ is the unique root in $\left(0,\,\frac{\pi \sigma^2}{b}\right)$ of the equation
	\[\frac{b\xi}{\sigma^2} = \arccos\left( \frac{\xi^2 + \mu(\mu+\theta \sigma^2)}{\sqrt{ (\xi^2 + \mu(\mu+\theta \sigma^2))^2 + \xi^2 \theta^2 \sigma^4}} \right).\]
	For $(\theta,\,\mu,\,b) \in \mathscr{B}_1$, $\psi(\theta) = 0$ and $h_\theta(x) \equiv 1$.
	Finally, for $(\theta,\,\mu,\,b) \in \mathscr{B}_2$, the solutions are given by $\psi(\theta) = -\frac{\mu^2}{2\sigma^2}$ and
	\[h_\theta(x) = e^{-\frac{\mu}{\sigma^2}x}\left( \frac{\mu}{\sigma^2}x +1 \right).\]
	The case of arbitrary $r_0$ and $r_b$ is conceptually similar, but requires even more complicated regions into which to separate the parameter space. For instance, it will be necessary
	to consider the signs of $\theta(r_0 + r_b)$, $(\mu-\theta r_0 \sigma^2)(\mu+\theta r_b \sigma^2)$, and $b(\mu-\theta r_0 \sigma^2)(\mu+\theta r_b \sigma^2) + \theta (r_0+r_b)\sigma^4$, amongst
	other quantities. It is therefore clear that an explicit description of the solution to \eqref{one-dim LD ODE} will, in general, be very complex.
  \end{example}

	\section{Acknowledgement}

  The first author gratefully acknowledges the mentorship and friendship of Professor Mikl\'{o}s Cs\"{o}rg\H{o}, over the years, and Professor Cs\"{o}rg\H{o}'s influence on both his
	research direction and academic career in the years that have passed since his graduation as a student at Carleton University.

  \clearpage
	\bibliographystyle{plainnat}
  \bibliography{myrefs}

  \end{document}